# Aitken delta-squared generalized Juncgk-type iterative procedure

M. De la Sen
Institute of Research and Development of Processes. University of Basque Country
Campus of Leioa (Bizkaia) – PO Box. 644- Bilbao, 48080- Bilbao. SPAIN
email: manuel.delasen@ehu.es

**Abstract**. This paper discusses a general Aitken delta-squared generalized Jungck-modified $S$ -iterative scheme. The study applies generalized versions of Aitken delta-squared procedure and Venter´s theorem to discuss positivity and global stability of the generalized Jungck iterative scheme which is of interest in numerical methods and its acceleration of convergence.

## 1. Introduction

Iterative schemes are of interest in numerical computations and its properties related to acceleration of convergence to solve scientific and engineering problems. A very common iterative scheme is the so – called Jungck iterative scheme, which involves the use of two coupled mappings, and its various extensions, [1-13]. Such an scheme is useful also in fixed point theory to find common fixed points of both mappings. In this paper, we extend such an iterative scheme and study its stability and positivity under certain parametrical restrictions. Also, we study the acceleration of convergence by proposing an Aitken type delta squared procedure for acceleration of convergence and combining the results with the generalization of a background Venter´s stability theorem result, whose basic form is well - known in discrete parametrical recursive identification , [14].

*1.1 Notation*

The real sequences $\{y_n\}$ and $\{x_n\}$ are equivalent, denoted by $\{y_n\} \approx \{x_n\}$, if both have the same limit.

## 2. The iterative sequence and some preliminaries

The following Aitken delta-squared generalized Jungck- modified $S$ -iterative scheme is of interest in order to accelerate the convergence of the modified sequence to the same limit as the unmodified one provided that such a limit exists:

$$Sz_{n+1} = (1-a_n)T^n z_n + a_n T^n y_n \ ; \ z_{-1} = 0 , \ \ z_0 \in C \tag{2.1.1}$$

$$ASz_{n+1} = Sz_{n+1} - \mu_{n+1} \frac{(\Delta S z_{n+1})^2}{\Delta^2 S z_{n+1}} \tag{2.1.2}$$

$$Sy_n = (1-b_n)Sz_n + b_n T^n z_n \tag{2.1.3}$$

$$ASy_n = Sy_n - \nu_n \frac{(\Delta S y_n)^2}{\Delta^2 S y_n} \tag{2.1.4}$$

; $\forall n \in \mathbf{N}$ , where $S, T : C \to C$ are two mappings on a nonempty subset $C$ of a Banach space $(X, \| \ \|)$ subject to $T(C) \subseteq S(C)$, $\{a_n\}$ and $\{b_n\}$ are real sequences in $[0 ,1]$, $\{\mu_n\}$ and $\{\nu_n\}$ are binary sequences

taking values 0 or 1; and $\Delta$ is a finite difference operator defining the correcting Aitken-type terms of the above iterative procedure as follows:

$$\Delta S z_{n+1} = S z_{n+2} - S z_{n+1} \;; \quad \Delta S y_n = S y_{n+1} - S y_n \tag{2.2.1}$$

$$\Delta^2 S x_{n+1} = S z_{n+1} - 2 S z_{n+2} + S z_{n+3} \;; \quad \Delta^2 S y_n = S y_n - 2 S y_{n+1} + S y_{n+2} \tag{2.2.2}$$

; $\forall n \in \mathbf{N}$ . The binary sequences $\{\mu_n\}$ and $\{\nu_n\}$ have two functions, namely, a) To remove the Aitken correction if it would imply division by zero at some iteration (i.e. if $\Delta^2 S z_{n+1}$ or $\Delta^2 S y_n$ is zero) ; b) To decide if implementing the Aitken correction or not at any particular iteration which does no imply division by zero.

The subsequent two technical preliminary results reflect the features that the sequences $\{S x_n\}$ and $\{S y_n\}$ can have different limits if such limits exist and that under mild conditions the Aitken correction leads to the same limits as its standard Jungck iterative process provided that such limits exist.

**Lemma 2.1**. The following properties hold:

**(i)** Assume that $\{a_n\} \subset [0, 1)$ and $\{S z_n\}$ converges to $(Sz)^*$ . Then, the following limit exists:

$$\lim_{n\to\infty} \left\{ [1 + a_n (b_n - 1)](Sz)^* - (1 - a_n) S y_n - a_n b_n T^n y_n \right\} = 0$$

**(ii)** Assume that $\{b_n\} \subset (0, 1]$ and $\{S y_n\}$ converges to $(Sy)^*$ . Then, the following limit exists:

$$\lim_{n\to\infty} \left\{ (1 - a_n)(Sz)^* - b_n S z_{n+1} - (1 - a_n)(1 - b_n) S z_n + a_n b_n T^n y_n \right\} = 0$$

**(iii)** Assume that $\{a_n\} \subset [0, 1)$, $\{b_n\} \subset (0, 1]$, $\{S z_n\}$ converges to $(Sz)^*$ and $\{S y_n\}$ converges to $(Sy)^*$ . Then, the following limit exists:

$$\lim_{n\to\infty} \left\{ [1 + a_n (b_n - 1)](Sz)^* - (1 - a_n)(S y)^* - a_n b_n T^n y_n \right\} = 0$$

If, in addition, $(Sx)^* = (Sy)^*$ then

$$\lim_{n\to\infty} \left\{ a_n b_n \left( (Sz)^* - T^n y_n \right) \right\} = 0$$

If $(Sz)^* = (Sy)^*$ and $\liminf_{n\to\infty} a_n b_n > 0$ then $\exists \lim_{n\to\infty} T^n y_n = (Sz)^*$ .

*Proof*: One gets obtaining explicitly $T^n z_n$ from (2.1.1), since $a_n \neq 1$, and after its replacement into (2.1.3), one gets:

$$b_n S z_{n+1} + (1 - a_n)(1 - b_n) S z_n = \left( (1 - a_n) S + a_n b_n T^n \right) y_n \;; \forall n \in \mathbf{N}$$

Property (i) follows since $\{S z_n\} \to (Sz)^*$ . To prove Property (ii), one obtains explicitly $T^n z_n$ from (2.1.3), since $b_n \neq 0$ , and one gets by replacing it into (2.1.1) that

$$(1 - a_n) S y_n - b_n S z_{n+1} - (1 - a_n)(1 - b_n) S z_n + a_n b_n T^n y_n = 0 \,, \ \forall n \in \mathbf{N}$$

Property (ii) follows since $\{Sy_n\} \to (Sy)^*$. Property (ii) follows from Properties (i) -(ii). □

**Lemma 2.2**. The subsequent properties follow:

**(i)** If $\{Sx_n\}$ converges to a finite limit $(Sz)^*$ then $\{ASz_n\} \to (Sz)^*$ provided that for any $n \in \boldsymbol{N}_0 = \boldsymbol{N} \cup \{0\}$

$$\Delta^2 Sz_{n+1} = (1-a_n)T^n z_n + a_n T^n y_n - 2(1-a_{n+1})T^{n+1}x_{n+1} - 2a_{n+1}T^{n+1}y_{n+1}$$
$$+ (1-a_{n+2})T^{n+2}x_{n+2} + a_{n+2}T^{n+2}y_{n+2} = 0 \Rightarrow \mu_{n+1} = 0$$

Furthermore, if $\{\mu_n\} \to 1$ then $\{ASz_n\} \approx \{Sz_n\}$ and $\{ASz_n\} \to (Sz)^*$ faster than $\{Sz_n\} \to (Sz)^*$.

**(ii)** If $\{Sy_n\}$ converges to a finite limit $(Sy)^*$ then $\{ASy_n\} \to (Sy)^*$ provided that for any $n \in \boldsymbol{N}$

$$\Delta^2 Sy_n = Sy_n - 2Sy_{n+1} + Sy_{n+2} = (1-b_n)Sz_n + b_n T^n x_n - 2(1-b_{n+1})Sx_{n+1} - 2b_{n+1}T^{n+1}x_{n+1}$$
$$+ (1-b_{n+2})Sz_{n+2} + b_{n+2}T^{n+2}x_{n+2} = 0 \Rightarrow \nu_n = 0$$

Furthermore, if $\{\nu_n\} \to 1$ then $\{ASy_n\} \approx \{Sy_n\}$ $\{ASy_n\} \to (Sy)^*$ faster than $\{Sy_n\} \to (Sy)^*$.

*Proof*: Note from direct computations that

$$\sum_{j=1}^{n}(A-I)(Sz_{j+1} - Sz_j) = \mu_{n+1}\frac{(Sx_{n+2} - Sz_{n+1})^2}{Sz_{n+1} - 2Sz_{n+2} + Sz_{n+3}} - \mu_1\frac{(\Delta Sz_1)^2}{\Delta^2 Sz_1}$$

; $\forall n \in \boldsymbol{N}$ and, since $\{Sy_n\}$ converges to a finite limit $(Sy)^*$, one gets:

$$\lim_{n\to\infty}\left[\sum_{j=1}^{n}(A-I)(Sz_{j+1} - Sz_j)\right] = -\mu_1\frac{(\Delta Sz_1)^2}{\Delta^2 Sz_1} = \lim_{n\to\infty}(A-I)(Sz_{n+1} - Sz_1)$$

$$= \lim_{n\to\infty}(A-I)Sz_{n+1} - (A-I)Sz_1 = A(Sz)^* - (Sz)^* - (A-I)Sz_1$$

$$= \lim_{n\to\infty} ASz_{n+1} - (Sz)^* - (A-I)Sz_1 = -\mu_1\frac{(\Delta Sz_1)^2}{\Delta^2 Sz_1}$$

$$= (A-I)Sz_1 = ASz_1 = -\mu_1\frac{(Sz_2)^2}{Sz_3 - 2Sz_2}$$

since $Sz_1 = 0$ from (2.1.1), $\Delta^2 Sz_{n+1} = 0 \Rightarrow \mu_{n+1} = 0$ and

$$(A-I)Sz_1 = -\mu_1\frac{(\Delta Sz_1)^2}{\Delta^2 Sz_1} = -\mu_1\frac{(Sz_2 - Sz_1)^2}{Sz_1 - 2Sz_2 + Sz_3} = -\mu_1\frac{(Sz_2)^2}{Sz_3 - 2Sz_2}$$

Then, $\{ASz_n\} \approx \{Sz_n\}$ and $\{ASz_n\} \to (ASz)^* = (Sz)^*$. On the other hand, simple calculations with (2.1.2) yield:

$$ASz_n - (Sz)^* = Sz_n - (Sz)^* - \mu_n\frac{(\Delta Sz_n)^2}{\Delta^2 Sz_n}$$

since $\left(ASz\right)^{*}=\left(Sz\right)^{*}$ so that if $\{\mu_n\}\to 1$ then $\{ASz_n\}\to \left(Sz\right)^{*}$ faster than $\{Sz_n\}\to \left(Sz\right)^{*}$ since $\left(\frac{(\Delta Sz_n)^2}{\Delta^2 Sz_n\left(Sz_n-\left(Sz\right)^{*}\right)}\right)\to 1$ as $n\to\infty$ so that the limit below exists:

$$\lim_{n\to\infty} \frac{ASz_n-\left(Sz\right)^{*}}{Sz_n-\left(Sz\right)^{*}}=1-\lim_{n\to\infty}\left(\frac{(\Delta Sz_n)^2}{\Delta^2 Sz_n\left(Sz_n-\left(Sz\right)^{*}\right)}\right)=0$$

Property (i) has been proved. Property (ii) is proved "mutatis-mutandis" from the evaluation of $\sum_{j=1}^{n}(A-I)\left(Sy_{j+1}-Sz_j\right)$; $\forall n\in \boldsymbol{N}$ so that

$$\lim_{n\to\infty}\left[\sum_{j=1}^{n}(A-I)\left(Sy_{j+1}-Sy_j\right)\right]=-\nu_1\frac{(\Delta Sy_1)^2}{\Delta^2 Sy_1}$$

$$=\lim_{n\to\infty} ASy_{n+1}-\left(Sy\right)^{*}-(A-I)Sy_1=(A-I)Sy_1=ASy_1=-\nu_1\frac{(Sy_2)^2}{Sy_3-2Sy_2} \qquad \square$$

The following result proves properties of stability and convergence for the case when $T,S:X\to X$ are linear and the sequences $\{a_n\}$ and $\{b_n\}$ are bounded while non-necessarily nonnegative.

**Lemma 2.3**. Assume that $\{a_n\}$ and $\{b_n\}$ are real bounded sequences and that $T,S:X\to X$ are linear and $S:X\to X$ is, furthermore, one-to-one and of closed range. Then, the following properties hold:

**(i)** If, for some $k_{1,2}\in \boldsymbol{R}_{0+}$, one has

$$\limsup_{n\to\infty}\left[\,|a_n|\left\|T^n\right\|-k_1\right]\le 0\,,\ \limsup_{n\to\infty}\left[\,|b_n|\left\|T^n\right\|-k_2\right]\le 0$$

$$\limsup_{n\to\infty}\left[\,|1-a_n|\left\|T^n\right\|-k_1'\right]\le 0\,,\ \limsup_{n\to\infty}\left[\,|1-b_n|-k_2'\right]\le 0$$

$$k_2'+\mu^{-1}(S)k_2\le 1\,;\ \mu^{-1}(S)\left[k_1'+k_1\left(k_2'+k_2\,\mu^{-1}(S)\right)\right]\le 1$$

with the inequality $\mu^{-1}(S)\left[k_1'+k_1\left(k_2'+k_2\,\mu^{-1}(S)\right)\right]\le 1$ being guaranteed if

$$k_2'+\mu^{-1}(S)k_2\le 1\,;\ \mu^{-1}(S)\left(k_1'+k_1\right)\le 1$$

Then, $\{z_n\}$ and $\{y_n\}$ are bounded. Furthermore, $\{Az_n\}$ and $\{Ay_n\}$ are also bounded under the conditions of Lemma 2.2.

**(ii)** If $\|T\|\le 1$ and

$$\mu^{-1}(S)\left(|1-b_n|\|S\|+|b_n|\right)\le M\,;\ \mu^{-1}(S)\left[|1-a_n|+\mu^{-1}(S)|a_n|\left(|1-b_n|\|S\|+|b_n|\right)\right]\le 1;\ \forall n\in \boldsymbol{N}$$

Then, $\{x_n\}$ and $\{y_n\}$ are bounded. Furthermore, $\{Az_n\}$ and $\{Ay_n\}$ are also bounded under the conditions of Lemma 2.2.

**(iii)** If $\|T\|\le 1$ and

$$\mu^{-1}(S)\left(|1-b_n|\|S\|+|b_n|\right)\le M\,;\ \left[|1-a_n|+|a_n|\,M\right]\le \mu(S);\ \forall n\in \boldsymbol{N}$$

Then, $\{x_n\}$ and $\{y_n\}$are bounded. Furthermore, $\{Az_n\}$ and $\{Ay_n\}$ are also bounded under the conditions of Lemma 2.2.

**(iv)** If $\|T\| < 1$ and $\{a_n\} \to 1$ then $\{z_n\}$and$\{y_n\}$ are bounded , equivalent and converge to zero. Under the conditions of Lemma 2.2, $\{Az_n\}$ and $\{Ay_n\}$ are bounded, equivalent and converge to zero faster than $\{z_n\}$and$\{y_n\}$.

**(v)** If $\|T\| < 1$ and $\{b_n\} \to 1$ then $\{z_n\} \approx \{y_n\}$ are bounded and converge to zero. Under the conditions of Lemma 2.2, $\{Ax_n\}$ and $\{Ay_n\}$ are bounded, equivalent and converge to zero faster than $\{x_n\}$and$\{y_n\}$.

*Proof*: The various properties for the Aitken-modified sequences$\{Az_n\}$ and$\{Ay_n\}$ which follow from the corresponding ones for $\{z_n\}$ and $\{y_n\}$ are direct under Lemma 2.2 so that no specific proofs are needed. Thus, only proofs for the generalized Jungck-modified $S$ -iterative scheme are now given. One gets from (2.1.1) and (2.1.3)

$$\|z_{n+1}\| \le \mu^{-1}(S)\Big( |1-a_n|\|T^n\|\|z_n\| + |a_n|\|T^n\|\|y_n\| \Big) \qquad (2.1.5)$$

$$\|y_n\| \le \mu^{-1}(S)\Big( |1-b_n|\|S\| + |b_n|\|T^n\| \Big)\|z_n\| \qquad (2.1.6)$$

where the minimum modulus of $S$ is positive, i.e. $\mu(S) > 0$ , since$S : X \to X$ is linear, one-to-one and of closed range then

$$\|z_{n+1}\| \le \mu^{-1}(S)\Big[ |1-a_n| + \mu^{-1}(S)|a_n|\Big( |1-b_n|\|S\| + |b_n|\|T^n\| \Big)\Big] \|T^n\|\|z_n\|$$

; $\forall n \in \boldsymbol{N}$ and one gets under the given assumptions and recursive calculations that

$$\|y_n\| \le \big(k_2' + \mu^{-1}(S)k_2\big)\|z_n\| \le \|z_n\| ;\ \forall n \in \boldsymbol{N}$$
$$\|z_{n+1}\| \le \mu^{-1}(S)\big[k_1' + k_1\big(k_2' + k_2\,\mu^{-1}(S)\big)\big]\ \|z_n\| \le \|z_n\| \le \|z_1\| ;\ \forall n \in \boldsymbol{N}$$

This proves the boundedness of Property (i) for the unmodified sequences and also for the Aitken delta-squared corrected sequences under the conditions of Lemma 2.2. Property (ii) follows from the subsequent relations which follow from the given assumptions:

$$\| z_{n+1}\| \le \|z_n\| \le \|z_1\| ;\ \|y_n\| \le \mu^{-1}(S)\big(|1-b_n|\|S\| + |b_n|\big)\|z_n\| \le M\|z_n\| \le M\|z_1\| ; \forall n \in \boldsymbol{N}$$

Property (iii) follows under sufficient conditions to prove the boundedness of the sequences under Property (ii) holds, by combining both given assumptions. Properties (iv)-(v) follow from (2.1.5), (2.1.6) with $\left\{\|T^n\|\right\} \to 0$ as follows:

$\{z_n\} \to 0$, $\limsup_{n\to\infty}\big[\|y_n\| - \mu^{-1}(S)\big(|1-b_n|\|S\|\big)\|z_n\|\big] \le 0$ so that$\{y_n\} \to 0$ if$\{a_n\} \to 1$; and $\{y_n\} \to 0$ and $\{z_n\} \to 0$ if $\{b_n\} \to 1$. □

## 3. Extended Venter´s Theorem and its use in the modified Juncgk iterative scheme

The subsequent preliminary result is then used:

**Theorem 3.1**. Consider the following iterative scheme parameterized by real sequences and constants:

$$x_{n+1} \le (1-\alpha_n+\gamma_n)x_n+\omega_n+\sigma \;;\; x_0 \ge 0 \tag{3.1}$$

where

$$\sigma \ge 0\,, \alpha_n \in (0\,,1]\,,\; \omega_n \ge 0\,,\; \{\gamma_n\} \text{ is bounded, } \gamma_n \ge 0\;;\; \forall n \in \boldsymbol{N}_0 = \boldsymbol{N} \cup \{0\} \tag{3.2.1}$$

$$\sum_{i=0}^{\infty}\omega_n < +\infty\,,\; \sum_{i=0}^{\infty}\alpha_n = +\infty\,,\; \alpha_n \to 0 \text{ as } n \to \infty \tag{3.2.2}$$

and define the real constant $K \in [0\,,1)$ as $K = \lim_{n\to\infty} \frac{1}{n+1}\sum_{i=0}^{\infty}(1-\alpha_n)$ (which exists since $\alpha_n \to 0$ as $n \to \infty$ ).Then, the following properties hold:

**(i)** If $\sigma = \gamma_n = 0$ ; $\forall n \in \boldsymbol{N}_0 = \boldsymbol{N} \cup \{0\}$ then $\{x_n\} \to 0$ (Venter´s theorem [* ] ). Furthermore,

$$\sum_{i=0}^{n} K^{n-i}\left[\omega_i - (K+\alpha_i-1)^i x_i\right] \ge -K^{n+1}x_0\;;\; \lim_{n\to\infty}\sum_{i=0}^{n} K^{n-i}\left[\omega_i - (K+\alpha_i-1)^i x_i\right] = 0 \tag{3.3}$$

**(ii)** If $\sigma = \gamma_n = 0$ ; $\forall n \in \boldsymbol{N}_0$ then $\{x_n\}$ is bounded , $\sum_{n=0}^{\infty}\alpha_n x_n < +\infty$ and $\sum_{n=0}^{\infty} x_n < +\infty$.

**(iii)** If $\sigma = 0$ ; $\forall n \in \boldsymbol{N}_0$ then $\{x_n\}$ is bounded and $\sum_{n=0}^{\infty}(\alpha_n - \gamma_n)x_n < +\infty$.

**(iv)** If $\limsup_{n\to\infty}(\gamma_n - \alpha_n) \le 0$ then $\{x_n\}$ is bounded. If, furthermore, $\sup_{n\in\boldsymbol{N}_0}(\gamma_n - \alpha_n) < 0$ then

$$\sup_{n\in\boldsymbol{N}_0} x_n \le \frac{1}{\inf_{n\in\boldsymbol{N}_0}(\alpha_n-\gamma_n)}\left((1-K)x_0 + \sigma + \sup_{0\le i\le n}\omega_i\right) < +\infty$$

*Proof*: If $\sigma = \gamma_n = 0$ ; $\forall n \in \boldsymbol{N}_0$ then $\{x_n\} \to 0$ from Venter´s theorem. The second part of Property (i) follows by rewriting (3.1) for $\sigma = \gamma_n = 0$ ; $\forall n \in \boldsymbol{N}_0$ in the equivalent form

$$x_{n+1} \le Kx_n + (1-\alpha_n-K)x_n + \omega_n + \sigma \;;\; x_0 \ge 0$$

leading to

$$x_{n+1} \le Kx_n + (1-\alpha_n-K)x_n + \omega_n \;;\; x_0 \ge 0$$

; $\forall n \in \boldsymbol{N}_0$ so that one gets via recursive calculations:

$$0 \le x_n \le K^n x_0 + \sum_{i=0}^{n-1}K^{n-i-1}\left[(1-\alpha_i-K)x_n + \omega_i\right]$$

with $K < 1$ which leads to (3.3) and Property (i) is proven. On the other hand, one gets directly from (3.1) for $\sigma = 0$

$$\lim_{n\to\infty}\sum_{i=0}^{n}(x_{i+1}-x_i) = -x_0 = \sum_{i=0}^{\infty}(\omega_i + (\gamma_i - \alpha_i)x_i)$$

which yields

$\sum_{n=0}^{\infty}(\alpha_n-\gamma_n)x_n=x_0+\sum_{n=0}^{\infty}\omega_n<+\infty$ and hence Property (iii). If, in addition, $\gamma_n=0$ ; $\forall n\in\boldsymbol{N}_0$ then $\sum_{n=0}^{\infty}\alpha_n x_n<+\infty$ , $0\le\left(\inf_{n\in\boldsymbol{N}_0}\alpha_n\right)\left(\sum_{n=0}^{\infty}x_n\right)\le\sum_{n=0}^{\infty}\alpha_n x_n<+\infty$ , $\sum_{n=0}^{\infty}x_n<+\infty$ and $\{x_n\}$ bounded. Hence, Property (ii). Note that, since $K<1$, (3.1) may be equivalently rewritten as follows:

$$
\begin{aligned}
x_{n+1} &\le Kx_n+(1+\gamma_n-\alpha_n-K)x_n+\omega_n+\sigma \\
&=K^{n+1}x_0+\sum_{i=0}^{n}K^{n-i}\left[(1+\gamma_i-\alpha_i-K)x_i+\omega_i+\sigma\right] \\
&\le K^{n+1}x_0+\frac{1-K^{n+1}}{1-K}\sup_{0\le i\le n}\left[(1+\gamma_i-\alpha_i-K)x_i+\omega_i+\sigma\right] \qquad (3.4.1)\\
&\le x_0+\frac{1-K^{n+1}}{1-K}\sup_{0\le i\le n+1}\left[(1+\gamma_i-\alpha_i-K)x_i+\omega_i+\sigma\right] \qquad (3.4.2)
\end{aligned}
$$

; $\forall n\in\boldsymbol{N}_0$ with $x_0\ge 0$ what implies

$$\sup_{0\le i\le n}x_i\le x_0+\sup_{0\le i\le n}\left[\frac{1-K^n}{1-K}\left((1+\gamma_i-\alpha_i-K)x_i+\omega_i+\sigma\right)\right] \qquad (3.5)$$

or

$$\left[1-K-\sup_{0\le i\le n}\left(1-K^i\right)\sup_{0\le i\le n}(1+\gamma_i-\alpha_i-K)\right]\sup_{0\le i\le n}x_i\le(1-K)x_0+\sigma+\sup_{0\le i\le n}\omega_i \qquad (3.6)$$

Then $\{x_n\}$ is bounded if $0\le\lim_{n\to\infty}\sup_{0\le i\le n}x_i<+\infty$ which is guaranteed from (3.6) if $\lim_{m>n}\sup_{n\to\infty}\left[\left(1-k^m\right)(\alpha_n-\gamma_n)\right]>0$ what is in turn guaranteed if $\lim\sup_{n\to\infty}(\gamma_n-\alpha_n)\le 0$ . Then, one gets from (3.4.1) that

$$\lim\sup_{n\to\infty}x_{n+1}\le\frac{1}{1-K}\lim\sup_{n\to\infty,}\left(\sup_{0\le i\le n}\left[(1+\gamma_i-\alpha_i-K)x_i+\omega_i+\sigma\right]\right)$$

If , furthermore, $\sup_{n\in\boldsymbol{N}_0}(\gamma_n-\alpha_n)<0$ then, it follows from (3.4.2) that

$$0\le\inf_{n\in\boldsymbol{N}_0}(\alpha_n-\gamma_n)\sup_{n\in\boldsymbol{N}_0}x_n\le(1-K)x_0+\sigma+\sup_{0\le i\le n}\omega_i<+\infty$$

with its first inequality being strict if $x_0>0$ . Thus, Property (iv) has been proven. □

**Theorem 3.2**. Consider the iterative scheme (2.1.1)-(2.14) under the subsequent constraints:

1) $\{b_n\}\subset(0,1]$, $\{b_n\}\to 1$

2) $\{a_n\}\subset[0,1]$, $\{a_n\}\to 1$, $\left\|T^n\right\|=o(|y_n|)$ or $\{a_n\}\subset[0,1]$, $\{a_n\}\to 0$, $\left\|T^n\right\|=o(|z_n|)$

3) $T,S:X\to X$ are linear and $S:X\to X$ is one-to-one with closed range with the sequence $\left\{S^{-1}T_n\right\}$consists of composite positive operators.

4) Define nonnegative real sequences $\{\gamma_n\}$ and $\{\alpha_n\}(\subseteq [0,1]) \to 0$, of general terms satisfying $\gamma_n = \frac{2-\alpha_n}{b_n}$, such that $\sum_{i=0}^{\infty} \alpha_i = \sum_{i=0}^{\infty} (2 - b_i \gamma_i) < \infty$.

Then, the iterative scheme (2.1.1)-(2.1.4) is globally stable for any $z_1 \geq 0, y_1 \geq 0$ while it fulfils $\{S z_n\} \to 0$ , and it is a nonnegative sequence , $\{S y_n\} \to 0$, $\{T^n y_n\} \to 0$, $\{T^n z_n\} \to 0$; and $\{AS z_n\} \to 0$ , respectively, $\{AS y_n\} \to 0$ faster than $\{S z_n\} \to 0$, respectively, $\{S y_n\} \to 0$.

*Proof*: The iterative scheme (2.1.1)–(2.1.4) verifies Theorem 3.1 with the definitions $S z_n \equiv x_n$ ; $n \in \boldsymbol{N}_{0+} = \boldsymbol{N} \cup \{0\}$ under the above constraints 1-4 with $\{\omega_n\}$ of general term

$$\omega_n = b_n^{-1} S y_n + a_n T^n (y_n - z_n) = b_n^{-1} (1 - b_n) S z_n + a_n T^n (y_n - z_n) + T^n z_n$$

$$= b_n^{-1} (1 - b_n) S z_n + a_n T^n S^{-1} (S y_n) + (1 - a_n) T^n S^{-1} (S z_n); \ \forall n \in \boldsymbol{N}_{0+}$$

being nonnegative for $z_1 \geq 0$ and $y_1 \geq 0$, summable, and convergent to zero. Then, the result follows from Theorem 3.1 using Lemmas 2.2-2.3. □

## ACKNOWLEDGEMENTS

The author is very grateful to the Spanish Government by its support of this research trough Grant DPI2012-30651, and to the Basque Government by its support of this research trough Grants IT378-10 and SAIOTEK S-PE12UN015. He is also grateful to the University of Basque Country by its financial support through Grant UFI 2011/07.

## REFERENCES

[1] I. Inchan, “Viscosity iteration method for generalized equilibrium points and fixed point problems of finite family of nonexpansive mappings”, *Applied Mathematics and Computation*, Vol. 219, Issue 6, pp. 2949-2959, 2012.

[2] M. De la Sen, “Stable iteration procedures in metric spaces which generalize a Picard-type iteration”, *Fixed Point Theory and Applications*, Vol. 2010, Article ID 572057, 15 pages, 2010.doi:10.1155/2010/953091.

[3] M. De la Sen, “Stability and convergence results based on fixed point theory for a generalized viscosity iterative scheme”, *Fixed Point Theory and Applications*, Article IDF 314581, Vol. 2009, doi: 10.1155/ 2009/314581, 2009.

[4] M. Delasen, “Online optimization of the free parameters in discrete adaptive –control systems”, *IEE Proceedings-D Control Theory and Applications*, Vol. 131, Issue 4, pp. 146-157, 1984.

[5] J.J. Minambres and M. De la Sen, “Application of numerical - methods to the acceleration of the convergence of the adaptive control algorithms- The one-dimensional case”, *Computers & Mathematics with Applications-Part A*, Vol. 12, Issue 10, pp. 1049-1056, 1986.

[6] M. Abbas, T. Nazir and S. Romaguera, “Fixed point results for generalizad cyclic contraction mappings in partial metric spaces”, *Revista de la Real Academia de Ciencias Exactas, Fisicas y Naturales*, *Serie A: Matematicas,*Vol. 106, No. 2, pp. 287-297, 2012.

[7] M. De la Sen, “Linking contractive self-mappings and cyclic Meir-Keeler contractions with Kannan self-mappings”, *Fixed Point Theory and Applications*, Vol. 2010, Article ID 572057, 23 pages, 2010.doi:10.1155/2010/572057.

[8] W.S. Du, “New cone fixed point theorems for nonlinear multivalued maps with their applications”, *Applied Mathematics Letters*, Vol. 24, No. 2, pp. 172-178, 2011.

[9] W. Laowang and B. Panyanak, “Common fixed points for some generalized multivalued nonexpansive mappings in uniformly convex metric spaces”, *Fixed Point Theory and Applications*, Article number 20, doi:10.1186/1687-1812-2011-20, 2011.

[10] H.K. Nashine and M.S. Khan, “An application of fixed point theorem to best approximation in locally convex space”, *Applied Mathematics Letters*, Vol. 23, Issue 2, pp. 121-127, 2010.

[11] A. Latif and M.A. Kutbi, “Fixed points for w-contractive multimaps”, *International Journal of Mathematics and Mathematical Sciences*, Vol. 2009, Article ID 769467, 7 pages, 2009, doi:10.1155/2009/769467.

[12] T. Husain and A. Latif, “Fixed points of multivalued nonexpansive maps”, *International Journal of Mathematics and Mathematical Sciences*, Vol. 14, No.3, pp. 421-430, 1991.

[13] M. S. Khan, “Common fixed point theorems for multivalued mappings”, *Pacific Journal of Mathematics*, Vol. 95, No. 2 , pp. 337-347, 1981.

[14] J.M. Mendel, Discrete Techniques of Parameter Estimation: The Equation Error Formulation, Marcel Dekker Inc., 1973.